\documentclass[11pt,draftcls,onecolumn]{IEEEtran}
\IEEEoverridecommandlockouts         

\usepackage{times} 
\usepackage{color}
\usepackage{graphicx}
\usepackage{setspace}
\usepackage{bbm}
\usepackage{mathdots,mathrsfs}
\usepackage{amssymb,latexsym,amsfonts,amsmath,cite}
\usepackage{stmaryrd}
\usepackage{caption}
\usepackage[dvips]{psfrag}
\newtheorem{theorem}{Theorem}
\newtheorem{lemma}{Lemma}

\newtheorem{definition}{Definition}
\newcommand{\R}{{\mathbb{R}}}

\newcommand{\G}{{\mathcal{G}}}

\newcommand{\NN}{\boldsymbol{\mathcal{N}}}

\newcommand{\Rel}{{\text{Re}}}
\newcommand{\diag}{{\text{diag}}}
\usepackage{pgf,tikz}

\DeclareMathOperator{\nul}{ker}

\begin{document}
\title{\huge  Systemic Measures for Performance and Robustness of Large--Scale Interconnected Dynamical Networks}


\author{Milad Siami$^\dag$ and Nader Motee$^\dag$ 
\thanks{This work is supported by the Office of Naval Research under award ONR N00014-13-1-0636.}
\thanks{$^\dag$ M. Siami and N. Motee are with the Department of Mechanical Engineering and Mechanics, Packard Laboratory, Lehigh University, Bethlehem, PA. Email addresses:  {\tt\small \{siami,motee\}@lehigh.edu} }}
\maketitle
\thispagestyle{empty}
\pagestyle{empty}

\begin{abstract}
In this paper, we develop a novel unified methodology for performance and robustness analysis of linear dynamical networks. We introduce the notion of systemic measures for the class of first--order linear consensus networks. We classify two important types of performance and robustness measures according to their functional properties: convex systemic measures and Schur--convex systemic measures. It is shown that a viable systemic measure should satisfy several fundamental properties such as homogeneity, monotonicity, convexity, and orthogonal invariance. In order to support our proposed unified framework, we verify functional properties of  several existing performance and robustness measures from the literature and show that they all belong to the class of systemic measures. Moreover, we introduce new classes of systemic measures based on (a version of) the well--known Riemann zeta function, input--output system norms, and etc. Then, it is shown that for a given linear dynamical network one can take several different strategies to optimize a given performance and robustness systemic measure via convex optimization. Finally, we characterized an interesting fundamental limit on the best achievable value of a given systemic measure after adding some certain number of new weighted edges to the underlying graph of the network. 
\end{abstract}

\section{Introduction}
\allowdisplaybreaks

The interest in control systems society for performance and robustness analysis of large--scale dynamical network is rapidly growing \cite{Bamieh12, Siami13cdc, Young10,Bamieh11, abbas, zelazo, LovisariGarinZampieriResistance, Nicola}. Improving global performance as well as robustness to external disturbances in large--scale dynamical networks are crucial for sustainability, from engineering infrastructures to living cells; examples include a group of autonomous vehicles in a formation, distributed emergency response systems, interconnected transportation
networks, energy and power networks, metabolic pathways and even financial networks. One of the fundamental problems in this area is to determine to what extent uncertain exogenous inputs can steer the trajectories of a dynamical network away from its working equilibrium point. To tackle this issue, the primary challenge is to introduce meaningful and viable performance and robustness measures that can capture essential characteristics of the network.  A proper measure should be able to encapsulate transient, steady--state, macroscopic, and microscopic features of the perturbed large-scale dynamical network.

In this paper, we propose a new unified methodology to classify proper performance and robustness measures for large--scale dynamical networks subject to external stochastic disturbance inputs. We take an axiomatic approach to quantify several essential properties of a sensible measure. We introduce the class of {\it systemic measures} and show that this class of measure should satisfy monotonicity, positive homogeneity, convexity, and orthogonal invariance conditions. It is shown that several existing and widely used  performance measures in the literature are in fact special cases of this class of systemic measures \cite{Siami14acc, zelazo,Young10,noami11,Jadbabaie13}.

The performance analysis of linear consensus  networks subject to external stochastic disturbances has been studied in \cite{Bamieh12, Siami13cdc,Young10,Siami13siam, SiamiNecSys, noami11, Jadbabaie13}, where the $\mathcal H_2$--norm of the network was employed as a scalar performance measure. In \cite{Bamieh12}, the authors interpret the $\mathcal H_2$--norm of the system as a macroscopic performance measure capturing the notion of coherence. It has been shown that if the Laplacian matrix of the underlying graph of the network is normal, the $\mathcal H_2$--norm is a function of the eigenvalues of the Laplacian matrix \cite{Young10}. In \cite{Siami13cdc}, the authors consider general linear dynamical networks and show that tight lower and upper bounds can be obtained for the $\mathcal H_2$--norm of the network from the exogenous disturbance input to a performance output, which are functions of the eigenvalues of the state matrix of the network. Besides the common $\mathcal H_2$--norm approach, there are several other performance measures that have been proposed in \cite{Bamieh12, zelazo, Olfati-saber07}.

{In this paper, we show that classes of system--norm, spectral, and entropy based performance and robustness measures enjoy similar functional properties.} These common properties enable us to identify and classify such measures under one umbrella, so called systemic measures. Particularly, we explore new connections between $\mathcal H_p$--norm (for range of exponents $1 \leq p \leq \infty$) of a first-order linear consensus network with (a version of) the well--known Riemann zeta function of the Laplacian matrix of the underlying graph of the network. We also characterize the class of entropy--based performance and robustness measures and show that this class of measures depends closely to the number of spanning trees in the underlying graph of the network.

The rest of the paper is organized as follows. Notation and basic notions are defined in Section II. The problem statement is addressed in Section III. By applying an axiomatic approach in Section IV, we characterize the class of systemic measures for consensus seeking networks. In Sections V and VI, some subclasses of systemic measures are studied. In Section VII, we focus on improving systemic performance and robustness of first-order linear consensus networks.


\section{Mathematical Preliminaries}
\allowdisplaybreaks

The sets of all positive and nonnegative real numbers are denoted by $\mathbb R _{++}$ and $\mathbb R _{+}$, respectively. Throughout this paper, it is assumed that all graphs are finite, simple, undirected and connected. A weighted graph $\G$ is represented by a triple $\G = (V(\G), E(\G),w)$, where $V(\G)$  is the set of nodes, $E(\G)\subset \big\{\{i,j\}~\big|~ i,j \in V(\G), ~i \neq j \big\}$ is the set of edges, and $w: E(\G) \rightarrow \R_{++}$ is the weight function. The degree of each node $i \in V(\G)$ is defined by 
\[d_i ~\triangleq ~\sum_{e=\{i,j\} \in E(\G)} w(e).\]
The adjacency matrix $ A = [a_{ij}]$ of graph $\mathcal{G}$ is defined by setting $a_{ij} = w(e)$ if $e=\{i,j\} \in E(\G)$, otherwise $a_{ij}=0$. 
The Laplacian matrix of $\mathcal G$ is defined by $L_{\G} \triangleq \Delta - A$, where $\Delta=\diag[d_{1},\ldots,d_{n}]$. The eigenvalues of $L_{\G}$ are indexed in ascending order $\lambda_1 \leq \lambda_2 \leq \cdots \leq \lambda_n $ and $\lambda_1=0$. The Moore-Penrose pseudo-inverse of $L_{\mathcal G}$ is denoted by $L_{\mathcal G}^{\dag}=[l_{ji}^{\dag}]$ which is a square, symmetric, doubly-centered and positive semidefinite matrix.

\begin{definition}
The {\it centering matrix} of size $n$ is defined by
\[M_n~\triangleq~I_{n} - \frac{1}{n}\mathbf 1_n \mathbf 1_n^T=I_{n} - \frac{1}{n}J_n,\]
where $I_{n}$ is the $n \times n$ identity matrix, $\mathbf 1_n$  the $n \times 1$ vector of all ones, and $J_n$  the $n \times n$ matrix of all ones.
\end{definition}

We denote the generalized matrix inequality with respect to the positive semidefinite cone by ``$\,\preceq \,$".
The beta function is defined by
\begin{eqnarray} 
\bold{\beta} (x,y)&=&\int_0^1 t^{x-1}(1-t)^{y-1}dt =\frac{\Gamma(x)\Gamma(y)}{\Gamma(x+y)}, 
\end{eqnarray}
where $\Rel \{x\}, \Rel \{y\} >0$, and $\Gamma(.)$ is the well--known gamma function.

\begin{definition}
A real--valued function $f$ is {\it permutation invariant} if and only if $f(x)=f(Px)$ for every permutation matrix $P \in \R^{n \times n}$.
\end{definition}

\begin{definition}
The real--valued function $f: \R_+^n \rightarrow \R$ is Schur--convex if $f(Dx) \leq f(x)$ for every doubly stochastic matrix $D$ and all  $x \in \mathbb R_+^n$.
\end{definition}

\section{Problem Statement}

We consider the class of first-order linear consensus networks over a weighted graph $\mathcal G$. Each node of the graph with index $i$ represents a subsystem with state variable $x_i \in \R$ for $i=1,\ldots,n$.  The state of the entire network is denoted by $x=\left[\begin{array}{cccc}x_{1} & x_{2} & \ldots & x_{n} \end{array}\right]^{T}$. Suppose that the dynamics of this class of networks are governed by 
\begin{eqnarray}
\NN(L_{\G};x_0):
\begin{cases}
\dot x(t)~=~-L_{\mathcal G} x(t)+\xi(t)\\
y(t)~=~M_{n}x(t)
\end{cases},
\label{first-order-G}
\end{eqnarray}
where $L_{\mathcal G}$ and $M_{n}$ are the Laplacian and centering matrices of $\G$, respectively. The exogenous disturbance input is denoted by $\xi(t)$ and the output of the network by $y(t)$. The exogenous disturbance input captures the effect of the uncertain environment on the dynamical network. One may represent  dynamical system \eqref{first-order-G} by symbol $\boldsymbol{\mathcal{N}}(L_{\G}; x_0)$, where $x_0$ is the initial condition of the network.  For a given fixed initial condition $x_0$, we can classify the set of all linear consensus networks $\boldsymbol{\mathcal{N}}(L_{\G}; x_0)$ and denote it by $\mathfrak{N}(x_0)$. Whenever it is not confusing, for simplicity of our notations we use notation $\boldsymbol{\mathcal{N}}(L_{\mathcal{G}})$ instead of $\boldsymbol{\mathcal{N}}(L_{\G}; x_0)$. 

The linear dynamical network \eqref{first-order-G} can be viewed as a system that has been already stabilized by a linear state feedback control law and operating in closed--loop. The sparsity pattern of the Laplacian matrix $L_{\G}$ is imposed by the topology of the underlying graph $\G$ and the corresponding weight function, which models the coupling structure and strength among the subsystems in the closed--loop system. The existence of such inherent sparsity--constraints on the topology of the underlying graphs play a foundational role in emergence of severe theoretical fundamental limits on the global performance and robustness of this class of dynamical networks. The impacts of such fundamental limits usually appear as
fundamental tradeoffs between various measures of performance and robustness in the presence of external disturbances, time--varying coupling structures,  and various modeling uncertainties.  Our main objective is to propose an unified approach to analyze performance and robustness of linear dynamical networks subject to stochastic exogenous disturbance inputs and quantify limits of performance and robustness due to the structure of the underlying graph of such networks.

\section{A Unified Framework for Systemic Measures}

We adopt an axiomatic approach to introduce and categorize a general class of performance and robustness measures that capture our intuition of a meaningful measure of performance and robustness in large--scale dynamical networks. Our approach characterizes several properties that a sensible performance and robustness measure should satisfy. Let us first define two basic algebraic operations on the space of linear consensus network.

\vspace{0.15cm}
\begin{definition}
For every given $\NN(L_{\G_1}), \NN(L_{\G_2}) \in \mathfrak{N}(x_0)$,  the addition and scalar multiplication operations on $\mathfrak{N}(x_0)$ is defined as follows: 
\begin{description}
\item[(i)] $\NN(L_{\mathcal{G}_1}) ~ + ~ \boldsymbol{\mathcal{N}}(L_{\mathcal{G}_2}) = \NN(L_{\mathcal{G}_1}+L_{\mathcal G_2})$
\item[(ii)] $\alpha \NN(L_{\G_1}) ~ = ~ \NN(\alpha L_{\G_1})$ for all positive scalars $\alpha$.
\end{description}
\end{definition}
\vspace{0.15cm}

The addition operation of two linear consensus networks is equivalent to the edge union operation on the underlying graphs of the two networks. The scalar  multiplication operation of a linear consensus network is equivalent to scaling the weight function of the underlying graph. 

\vspace{0.15cm}
\begin{definition}[Convex Systemic Measures]
\label{def-systemic}
For a given space of linear networks $\mathfrak{N}(x_0)$, a convex systemic measure is an operator $\rho: \mathfrak{N}(x_0) \rightarrow \R$ with the following  properties for all~ $\NN (L_{\G _1}), \NN (L_{\G _2}) \in \mathfrak{N}(x_0)$: 
\begin{description}
\item[(i)] {\bf Positive homogeneity of degree $-1$:} for all ~$\kappa >0$ 
\[\rho \big (\kappa \NN(L_{\G_1}) \big )~=~{\kappa}^{-1} \rho \big (\NN(L_{\G_1})\big ),\]
\item[(ii)] {\bf Monotonicity:} ~If ~$L^{\dag}_{\mathcal G_1} \preceq L^{\dag}_{\mathcal G_2}$ ~then ~
\[\rho \big (\NN (L_{\G _1})\big ) ~\leq~ \rho \big (\NN (L_{\G_2}) \big),\]

\item[(iii)] {\bf Convexity:}~for all ~$0 \leq \alpha \leq 1$
\begin{eqnarray*}
& & \hspace{-0.7cm} \rho \big (\NN(\alpha L_{\mathcal G_1}+(1-\alpha)L_{\mathcal G_2}) \big )~\leq~  \\
& & \hspace{1.2cm} {\alpha}\rho \big (\NN(L_{\G_1}) \big )+{(1-\alpha)}\rho \big (\NN(L_{\G_2})\big ).
\end{eqnarray*}
\end{description}
\end{definition}
\vspace{0.15cm}

The monotonicity property imposes a partial ordering on the space of networks $\mathfrak{N}(x_0)$ (see \cite{Siami14acc} for some related discussions) and implies that a systemic measure $\rho$ is {\it subadditive} over the set of all linear consensus networks, i.e., 
\[\rho \big (\NN(L_{\G_1}) + \NN(L_{\G_2}) \big )\leq \rho\big ( \NN(L_{\G_1})\big )+\rho\big (\NN(L_{\G_2})\big ),\]
for all~ $\NN (L_{\G _1}), \NN (L_{\G _2}) \in \mathfrak{N}(x_0)$. This property can be interpreted as a fundamental tradeoff between systemic measures and sparsity of the underlying graph. If we add more edges to an existing graph, the value of the systemic measure will decrease. For instance in Theorem 2 of \cite{Siami14arxiv}, we explicitly show this relationship between the $\mathcal{H}_2$--norm of the system from $\xi$ to $y$ and the sparsity of the underlying graph. The homogeneity property implies that among all graphs with identical interconnection topologies, the ones with larger (stronger) coupling weights have smaller systemic measures. 

In some applications in dynamical networks, the desired performance measures may not be positively homogeneous of degree $-1$. In these situations, we can relax Definition \ref{def-systemic} by removing the homogeneity property and replacing it by an orthogonal invariance property. 

\vspace{0.15cm}
\begin{definition}[Schur--Convex Systemic Measures]
\label{def-schur-systemic}
For a given space of linear networks $\mathfrak{N}(x_0)$, a schur--convex systemic measure is an operator $\rho: \mathfrak{N}(x_0) \rightarrow \R$ that satisfies properties (ii) and (iii) in Definition  \ref{def-systemic} and is orthogonally invariant, i.e., 
\[ \rho \big( \NN(L_{\G}) \big)~=~\rho \big( \NN(UL_{\G}U^T) \big), \]
for all orthogonal matrices $UU^T=U^TU=I_n$.
\end{definition}
\vspace{0.15cm}
A Schur--convex systemic measure is a permutation invariant function of the Laplacian eigenvalues \cite{Siami14acc}. If a real--valued function is convex and permutation invariant, then it is a Schur--convex function  \cite{marshall11}. This implies that all orthogonally invariant convex systemic measures are Schur--convex systemic measures \cite{Siami14acc}, but vice versa is not always true. 

Some important examples of convex and Schur--convex systemic measures are summarized in Table \ref{table}. In the following sections, we will classify general classes of such systemic measures and show that our unified framework provides convex and tractable formulations  to optimize systemic measures for the class of linear consensus networks. 

\begin{table*}[t]
\centering
\small{
\begin{tabular}{ c | c | c  }  
\hline
${\rho}(.)$& Schur--convex systemic measure & Convex systemic measure  \\
\hline
\hline
Convergence time of the first--order&\\ consensus networks: $\frac{1}{\lambda_2}$& \checkmark & \checkmark \\
\hline
Laplacian energy of the first--order&\\ consensus networks: $\sum_{i=2}^n\frac{1}{2\lambda_i}$& \checkmark & \checkmark \\
\hline
Laplacian energy of the second--order&\\ consensus networks: $\sum_{i=2}^n\frac{1}{2\lambda_i^2}$& \checkmark & \\
\hline
Normalized higher order Laplacian energy of &\\ consensus networks: $\Big (\sum_{i=2}^n\frac{1}{\lambda_i^p}\Big)^{\frac{1}{p}}$& \checkmark & \checkmark\\
\hline
Local error of first-order&\\ consensus dynamics: $\frac{1}{2}\sum_{i \in V(\G)} \frac{1}{d_i}$&& \checkmark \\
\hline
$\mathcal H_p$--norms of first-order&\\ consensus networks: for all $1 \leq p \leq \infty$&\checkmark&\\
\hline
Entropy of the first--order&\\ consensus networks: $-\sum_{i=2}^n\log \lambda_i$& \checkmark &\\
\hline
\end{tabular}}
\caption{\small{Examples of convex systemic measures and Schur--convex systemic measures. }}
\label{table}
\end{table*}

\section{Convex Systemic Measures}
\label{section-spectral}

Our focus will be on two important classes of convex systemic measures. First, we investigate convex performance and robustness measures that are defined using spectral properties of the underlying graph. Next, we consider a class of convex systemic measure that is defined based on spatial specifications of the underlying graph. In the following, we discuss that these two seemingly different classes of measures enjoy similar fundamental properties as described in Definition \ref{def-systemic}. 
\subsection{Spectral--Based Systemic Measures}\label{subsect-spectral}
In this subsection, we classify an important class of convex systemic measures that are defined using Laplacian eigenvalues of the underlying graphs. Several well--known and widely used performance and robustness measures for linear consensus networks are indeed convex systemic measures (see Table \ref{table}). In the following, we identify a general subclass of orthogonally invariant convex systemic measures based on spectral zeta function.
\vspace{0.15cm}
\begin{definition}\label{Zeta}
For a given Laplacian matrix $L_{\G}$, the corresponding {\it spectral zeta function} is a complex--valued function  and is defined by
\[\boldsymbol {\zeta}_{\G}(p)~\triangleq~\sum_{\lambda_i \neq 0} \lambda_i^{-p},\]
where $\lambda_i$'s are  eigenvalues of the Laplacian matrix and $p$ is a real number.
\end{definition}

\vspace{0.15cm}
\begin{theorem}
\label{th-zeta1}
For some given parameters $1 \leq {p} \leq \infty$ and $k > 0$, the following spectral--based measure  
\begin{equation}
\rho\big(\NN(L_{\G})\big)\,=\,k \left(\boldsymbol{\zeta}_{\G}(p)\right)^{\frac{1}{p}}, \label{zeta-measure}
\end{equation}
is {an orthogonally invariant} convex systemic measure  for all $\NN(L_{\G}) \in \mathfrak{N}(x_{0})$. 
\end{theorem}

The spectral--based systemic measure (\ref{zeta-measure}) includes several well--known performance and robustness measures as its special cases. We discuss some of these cases. When  $p \rightarrow \infty$, we have 
\begin{eqnarray}
\rho\big(\NN(L_{\G})\big)&=&k \lim_{p  \rightarrow \infty}\left(\boldsymbol{\zeta}_{\G}(p)\right)^{\frac{1}{p}}~=~ \frac{k}{\lambda_2}.\nonumber 
\end{eqnarray}
In this case, our proposed spectral--based systemic measure reduces to the rate of convergence of the consensus process in dynamical network $\NN(L_{\G})$. On the other hand, if  $p \rightarrow -\infty$, our systemic measure boils down to 
\begin{eqnarray}
\rho\big(\NN(L_{\G})\big)&=&k \lim_{p  \rightarrow -\infty}\left(\boldsymbol{\zeta}_{\G}(p)\right)^{\frac{1}{p}}~=~ \frac{k}{\lambda_n}.\nonumber
\end{eqnarray}
For $p =1$, our proposed spectral--based systemic measure is exactly equal to the first--order Laplacian energy of $\NN(L_{\G})$ with exogenous white Gaussian noise input with identity covariance \cite{Siami14arxiv}. In this case, we have
\begin{eqnarray}
\rho\big(\NN(L_{\G})\big) 
& = & \sum_{i=2}^n \frac{1}{2\lambda_i}~= ~\frac{1}{2} \boldsymbol{\zeta}_{\G}(1), \nonumber
\end{eqnarray}
which is indeed equal to the $\mathcal{H}_2$--norm of the system from the exogenous disturbance input to the output. When $p = 2$, our proposed systemic measure is 
equal to the second--order Laplacian energy of a second-order linear consensus network. For an extensive discussion on this case, we refer the reader to \cite{Bamieh11,Siami13cdc, Siami14arxiv}.

\section{Schur--Convex Systemic Measures}\label{subsection-h-p}

%
In this section, we turn our attention to two important classes of Schur--convex systemic measures.  Some important examples of such Schur--convex systemic measures are summarized in Table \ref{table}. 

\subsection{$\mathcal H_p$--Based Systemic Measures}\label{subsec-H-p}

This class of systemic measures is defined using the Schatten $p$--norm of a matrix \cite{hornJohnson90}
\[\|A\|_{p^{\ast}}~=~\left(\sum_{i=1}^n \sigma_i^p \right)^{\frac{1}{p}},\]
where $\sigma_i$'s  are singular values of $A$ and $1 \leq p \leq \infty$. The Schatten $p$--norms are unitary invariant norms. When $p = 2$, the Schatten norm reduces to the well--known Frobenius norm of a matrix. For $p =\infty$, the Schatten norm is equivalent to the spectral norm, i.e., the induced $2$--norm of a matrix.

\vspace{0.1cm}
\begin{theorem}
For a given linear consensus network $\NN(L_{\G}) \in \mathfrak{N}(x_0)$, let us define the input-output $\mathcal H_p$--norm of the network for every $1 \leq p \leq \infty$  by
\[ \|G(s)\|_{\mathcal H_p}~=~ \left( \frac{1}{2\pi} \int_{-\infty}^{\infty}\| G(j\omega)\|_{{p^{\ast}}}^p \hspace{0.05cm} d\omega \right )^{\frac{1}{p}}\]
where $G(s)$ is the transfer function of the network  (\ref{first-order-G}) from $\xi(t)$ to $y(t)$ and $\|\hspace{0.03cm}.\hspace{0.03cm}\|_{p^{\ast}}$ is the Schatten $p$-norms. Then, 
\begin{equation}
\rho \big (\NN(L_{\G})\big )~=~\|G(s)\|_{\mathcal H_p} \label{sys-meas-H-p}
\end{equation}
 is a Schur--convex systemic measure for all $1 \leq p \leq \infty$. 
\end{theorem}
The $\mathcal H_p$--norm based systemic measures captures several important performance and robustness features of large--scale dynamical networks. For exponent $p=2$, the systemic measure \eqref{sys-meas-H-p} is equivalent to the input-output $\mathcal H_{2}$--norm of the network
\begin{equation}
 \|G(s)\|_{\mathcal H_{2}}~=~\Big(\sum_{i=2}^n\frac{1}{2\lambda_i}\Big)^{\frac{1}{2}},
\end{equation}
where $\lambda_i$'s are the eigenvalues of $L_{\G}$. This systemic measure quantifies to what extend the effect of exogenous stochastic disturbance inputs  propagate throughout the network \cite{Siami14arxiv,zelazo}. It can also capture a notion of coherence in linear consensus  networks \cite{Bamieh12}.  At the other end of the spectrum when $p=\infty$, the systemic measure \eqref{sys-meas-H-p} is equivalent to the  input-output $\mathcal H_{\infty}$--norm of the network
\begin{equation}
 \|G(s)\|_{\mathcal H_{\infty}}~=~\frac{1}{\lambda_2},
\end{equation}
where $\lambda_2$ is the second largest eigenvalue of $L_{\G}$, i.e., the algebraic connectivity of graph $\G$. This measure can be viewed as the maximum system gain when inputs are taken over all measurable signals with finite energy, i.e., input signals in $L^2([0,\infty);\R^n)$. In this case, the corresponding systemic measure carries important information about the worst--case input that can deteriorate the performance of the network significantly. Moreover, this systemic measure has implications for disturbance rejection and can be viewed as a measure of robust stability. The following result shows that there is a close relationship between $\mathcal H_p$--norm based systemic measures and the spectral--based convex systemic measures.

\begin{theorem}
\label{th-zeta}
For a given linear consensus network $\NN(L_{\G}) \in \mathfrak{N}(x_0)$, we have 
\begin{equation}
 \|G(s)\|_{\mathcal H_p}^p~=~\frac{-1}{{\beta} (\frac{p}{2},-\frac{1}{2})} \hspace{.05cm} \boldsymbol {\zeta}_{\G}(p-1)
 \end{equation}
for every $1 \leq p \leq \infty$, where  $\boldsymbol {\zeta}_{\G}(.)$ is the zeta function of the underlying graph $\G$ and ${\beta} (.,.)$ is the Beta function.
\end{theorem}

The result of  Theorem \ref{th-zeta} asserts that there is an inherent relationship between the system-theoretic properties large--scale dynamical networks and the structural properties of the underlying graph of the network. The zeta--function of a graph can be related to various characteristics of the graph and it can be shown that how it scales with the network size. We refer the reader to \cite{Siami14arxiv} for an extensive discussion. 
\subsection{Systemic Measures Generated by Schur--Convex Sums}
The second important class of Schur--convex systemic measures is generated by sums of convex decreasing functions (also known as Schur--convex sums) of Laplacian eigenvalues. 

\begin{theorem}
\label{f-sum}
Suppose that $f: \R_+ \rightarrow \R$ is a decreasing convex function. For every $\NN(L_{\G}) \in \mathfrak{N}(x_{0})$, the class of measures that are defined by
\begin{equation}
\rho\big( \NN(L_{\G})\big)~=~\sum_{i=2}^n f(\lambda_i),
\end{equation} 
are Schur--convex systemic measure. 
\end{theorem} 
%
Several examples of well--known performance and robustness measures that belong to the class of Schur--convex systemic measures are listed in Table \ref{table}. The first-- and second--order Laplacian energies are studied in detailed in \cite{Siami14arxiv}. In order to show how a systemic measure can be related to the structural properties of the underlying graph of the network, we focus on the following Schur--convex systemic measure 
\begin{equation}
\rho \big(\NN(L_{\G})\big)~=~ - \sum_{i=2}^n \log \lambda_i.
\label{entropy-measure}
\end{equation}
This measure is also known as an {\it entropy} measure for linear consensus networks \cite{Siami14acc}.  The systemic measure (\ref{entropy-measure}) can be interpreted as the logarithm of minimum--volume ellipsoid covering  the projection of steady--state output vectors of $\NN(L_{\G})$ along $\nul(L_{\G})$ onto $\nul(L_{\G})^{\perp}$.  
Let us denote the total number of spanning trees of  the underlying graph $\G = (V(\G), E(\G),w)$ of the network by 
\[ \boldsymbol{\tau}({\mathcal G})~ \triangleq~ \sum_{\mathcal T}\prod_{e \in E(\mathcal T)} w(e), \]
where the summation runs over all spanning trees $\mathcal T$ of $\mathcal G$.

\begin{lemma}
For a given linear consensus network $\NN(L_{\G}) \in \mathfrak{N}(x_0)$, the systemic entropy measure (\ref{entropy-measure}) can be calculated by 
\begin{eqnarray}
\rho \big(\NN(L_{\G})\big)~=~\log \left(\frac{n}{\boldsymbol{\tau}(\mathcal G)}\right),
\label{tree-eq}
\end{eqnarray}
where $\boldsymbol{\tau}(\mathcal G)$ is the total number of spanning trees of  $\G$ and $n$ is the number of nodes.
\end{lemma}

\section{Convex Optimization Based Formulations to Improve Systemic Measures}

In this section, we formulate several convex optimization problems in order to design network topologies with minimal systemic measures. Specifically, we consider three interesting scenarios for minimization of systemic measures by: adjusting edge weights in a dynamical network with fixed topology, rewiring the underlying graph, and adding new edges to the underlying graph of the network. In the following subsections, we will discuss  these cases in more details.  
\subsection{Adjusting Edge Weights in Dynamical Networks with Fixed Topologies.}
We investigate the problem of allocating new additional weights to some edges of a graph of a network in order to minimize a given convex or Schur--convex systemic measure subject to the constraint that the sum of all allocated weights have to add up to a given constant. This  constant can be normalized to number one. It is known \cite{boyd2006} that when $\rho$ is a permutation invariant closed convex function of Laplacian eigenvalues, then ${\rho}$ can be rewrite as a convex function of edge weights $\mathbf w=[w_{1} ~ w_{2} ~ \ldots ~ w_{m} ]^{T}.$ This implies that for Schur--convex systemic measures our design problem can be cast as a convex optimization problem and solved efficiently in polynomial time. Suppose that we are given a systemic measure $\rho$ that is defined as a real--valued function of the pseudo-inverse of the Laplacian matrix of the underlying graph. The problem of adjusting edge weights in a dynamical network with fixed topology can be cast as 
\begin{eqnarray}
\begin{array}{lcc}
\text{\bf Minimize} \quad {\rho}\big( \NN(L_{\mathcal G}) \big)={\phi}( L_{\G}^{\dag} )\nonumber \\
\text{\bf subject to:} \quad \mathbf w \geq 0, \quad \mathbf 1_m^{T} \mathbf w=1
 \end{array}.
\end{eqnarray}

Let us now consider the following auxiliary Semidefinite programing problem (see \cite{Ghosh:2008:MER} for similar techniques):
\begin{equation} \label{sdp}
\begin{array}{lcc}
\text{\bf Minimize} \quad  {\rho}\big( \NN(L_{\mathcal G}) \big)= \phi(Y-\frac{1}{n}J_n)\nonumber \\
\text{\bf subject to:} \quad \mathbf 1_m^{T} \mathbf w=1, \quad \mathbf w \geq 0,\nonumber\\
~~~~~~~~~~~~~\left[ \begin{array}{ccc} L_{\G}+\frac{1}{n}J_n & I_n  \\ I_n & Y  \end{array} \right] \succcurlyeq 0,
 \end{array}
 \end{equation}
where $\mathbf{w} \in \R^m$ and $Y\in \R^{n \times n}$ is the slack symmetric matrix. In order to show that these two problems are equivalent, we look at the Schur complement of block matrix
\begin{equation*}
\left[ \begin{array}{ccc} L_{\G}+\frac{1}{n}J_n & I_n  \\ I_n & Y  \end{array} \right] \succcurlyeq 0, 
\end{equation*}
which is equivalent to $Y-{1}/{n}J_n\succcurlyeq L_{\G}^{\dag}$.
According to the monotonicity property of a systemic measure, one can conclude that minimizing the Semidefinite programming problem minimizes the original problem with cost function $\phi( L_{\G}^{\dag})$. 

There are several important Schur--convex systemic measures that can be written as a function of the pseudo-inverse of the Laplacian matrix of the underlying graph. For instance, let us consider the problem of minimizing the first--order Laplacian energy of a linear consensus  network by adjusting the edge weights while the topology of the underlying graph is fixed (see \cite{Siami14acc} for more details). For a total effective resistance interpretation of this minimization problem, we refer the reader to \cite{Ghosh:2008:MER}. We can cast this design problem as
\begin{equation}
\begin{array}{lcc}
\text{\bf Minimize} \quad {\rho}\big( \NN(L_{\mathcal G}) \big)=\frac{1}{2}\mathbf{Tr} L_{\G}^{\dag} \nonumber \\
\text{\bf subject to:} \quad \mathbf w \geq 0, \quad \mathbf 1_m^{T} \mathbf w=1.
 \end{array}
\end{equation}
With our proposed reformulation technique, we can equivalently solve the following convex optimization problem to minimize the Laplacian energy
\begin{equation}
\begin{array}{lcc}
\text{\bf Minimize} \quad {\rho}\Big( \NN(L_{\mathcal G}) \Big)=\frac{1}{2}\mathbf{Tr} Y - \frac{1}{2} \nonumber \\
\text{\bf subject to:} \quad \mathbf 1_m^{T} \mathbf w=1, \quad \mathbf w \geq 0,\nonumber\\
~~~~~~~~~~~~~\left[ \begin{array}{ccc} L_{\G}+\frac{1}{n}J_n & I_n  \\ I_n & Y  \end{array} \right] \succcurlyeq 0,
 \end{array}
\end{equation}
with $w \in \R^m$ and $Y\in \R^{n \times n}$ is the slack symmetric matrix.

\subsection{Rewiring the Underlying Graph of the Network} 
In the second scenario, we focus on the problem of rewiring the underlying graph of a linear consensus network in order to minimize a given Schur--convex systemic measure. The total number of edges and their weights that can participate in rewiring is given and fixed.    
Suppose that the set of all simple connected graphs with $n$ nodes, $m$ edges and $\sum_{e \in E(\mathcal G)} \omega(e) = \alpha$ by $\mathbb M_{n,m,\alpha}$. Our network design problem can be cast as 
\begin{eqnarray}
\begin{array}{lcc}
\text{\bf Minimize} \quad {\rho}\big(\boldsymbol{\mathcal N}(L_{\mathcal G})\big)\nonumber \\
\text{\bf subject to:} \quad \mathcal G \in \mathbb M_{n,m,\alpha}
\end{array},
\end{eqnarray}
where $\rho$ is a given Schur--convex systemic measure. It can be shown that when all edge weights are equal, the  resulting graph after removing $k$ disjoint edges from $\mathcal K_n$ (complete graph with $n$ nodes) minimizes all Schur--convex measures over all simple connected graphs with $n$ nodes and $\frac{n(n-1)}{2}-k$ edges. For example, among all linear consensus networks with unweigted underlying graphs with $4$ nodes and $4$ edges, those networks with cyclic topologies minimize all Schur--convex systemic measures.
\subsection{Adding New Edges to the  Graph of the Network}
We limit our discussion only to linear consensus networks $\boldsymbol{\mathcal{N}}(L_{\mathcal G}) \in \mathfrak{N}(x_{0})$ endowed with the following class of  Schur--convex systemic measures
\begin{equation}
\rho \big(\NN(L_{\mathcal G})\big)~=~\sum_{i=2}^n f(\lambda_i),
\label{meas}
\end{equation}
where $f$ is a real--valued decreasing convex function and $\lim_{x \rightarrow \infty} f(x) = 0$. We refer to Table \ref{table} for some examples of this class of systemic measures. 
Let us denote the resulting network after adding at most $k$ edges to the underlying graph of the network by $\NN \big( L_{\G}^{new}\big)$. The following result characterizes a {\it fundamental limit} on the best achievable value for the systemic measure after adding at most $k$ new arbitrary weighted edges to the underlying graph of the network. 
\vspace{0.15cm}
\begin{theorem}\label{w}
Suppose that linear consensus network $\boldsymbol{\mathcal{N}}(L_{\mathcal G}) \in \mathfrak{N}(x_{0})$ is endowed with the performance and robustness systemic measure  \eqref{meas}. There is a fundamental limit on the best achievable performance and robustness systemic measure through adding at most $k$ new arbitrary weighted edges to the graph $\mathcal G$ of the network, i.e., 
\begin{equation}
\rho \big( \NN ( L_{\G}^{new}) \big) \geq \sum_{i=k+2}^{n}f(\lambda_i),
\label{result}
\end{equation}
 where $\lambda_i$'s are the Laplacian eigenvalues of the original graph of the network before adding new edges.  
\end{theorem}
%
\section{Conclusion}
In this paper, we proposed a new unified mathematical framework to study performance and robustness analysis in linear dynamical networks. Our main focus was on the class of first--order linear consensus networks. We introduced the notion of systemic measures for this class of networks. Two major classes of performance and robustness measures were classified: convex systemic measures and Schur--convex systemic measures. Depending on the application, it was discussed that a viable systemic measure should satisfy several fundamental properties such as homogeneity, monotonicity, convexity, and orthogonal invariance. Then, we showed that for a given linear dynamical network one can take several different strategies to optimize a given performance and robustness systemic measure. Several convex optimization problems were formulated to minimize a given systemic measure by: adjusting edge weights in a dynamical network with fixed topology, rewiring the underlying graph, and adding new edges to the underlying graph of the network. Finally, we characterized an interesting fundamental limit on the best achievable value of a given systemic measure after adding some curtain number of new weighted edges to the underlying graph of the network.

\begin{spacing}{.887}
\bibliography{references-Milad-7-2014} 
\end{spacing}
\end{document}